\documentclass{mcom-l}

\usepackage{amssymb}
\usepackage{graphicx}
\usepackage{xcolor}

\makeatletter
\@namedef{subjclassname@2020}{Subject}
\makeatother

\copyrightinfo{2024}{Eric Easthope}

\theoremstyle{definition}

\theoremstyle{remark}

\numberwithin{equation}{section}

\begin{document}
\emergencystretch 3em
\title[Fregean Flows]{Fregean Flows}

% Only \author and \address required
\author[Eric Easthope]{Eric Easthope}
\address{University of British Columbia, Vancouver, British Columbia, Canada}
\email{mail@eric.cyou}

% \subjclass required
\date{}
\dedicatory{This paper is dedicated to its first adversary Prof. Dave Gilbert.}

\begin{abstract}
I humbly introduce a concept I call \textit{Fregean flows}, a graph theoretic representation of classical logic, to show how higher-dimensional graph characteristics might be useful to prove or perhaps at best show the provability of simple deductive statements typically represented as one-dimensional strings of characters. I apply these to a very simple proof, namely proving the equivalence of two definitions for an Abelian group $G$, an ``if-and-only-if'' statement, using a re-representation of statements as vertices and both conjunctions and implications as differently ``coloured'' edges. This re-representation of an if-and-only-if is simple but shows unexpected geometry, and I discuss its possible utility in terms of provability through ideas of graph topology, similarities of graph contraction to deductive elimination, and recursion.
\end{abstract}

\maketitle

\begin{quotation}
\small
``If the task of philosophy is to break the domination of words over the human mind [...] then my concept notation, being developed for these purposes, can be a useful instrument for philosophers [....] I believe the cause of logic has been advanced already by the invention of this concept notation.'' \textemdash\, Gottlob Frege (preface to \textit{Begriffsschrift}, 1879; English translation)
\end{quotation}

\pagebreak

\section{Motivation}

\subsection{Logic as picture}

Gottlob Frege's visual \textit{Begriffsschrift} (1879) \cite{Frege} arguably canonized the foundations of first-order axiomatic logic and correspondingly the foundations for mathematical proof by symbolic logical reasoning as well. These early writings in analytical philosophy were highly visual, recursive, and not so much adhered to the one-dimensional symbol-sentence style of modern proofs. Frege was keen to develop purely logical principles of inference that could represent mathematical proof axiomatically without gaps and unwritten appeals to ``intuition,'' and he did so without sentences, higher-dimensionally, with Boolos describing Frege's ``two-dimensional symbolism [as depicting] the logical structure of statements visually and vividly.'' \cite{Boolos}. It is unfortunate that Boolos went on to call Frege's work ``nothing so much as wallpaper'' \cite{Boolos} without elaborating on what he indicated he knew about Frege's interest to develop arithmetic through geometry, but Boolos' criticism highlighted an emergent and rather novel adherence to logic represented as sentence over logic represented as picture.

\subsection{The non-linearity of proof}

Text, and correspondingly the tasks of reading and writing text, are one-dimensional. There really is nothing that precludes all written word from being expressed as an endless, one-dimensional string of characters. And this how we express threads of logic too, not to mention threads of executed machine code. Pragmatically we represent text linearly to be maximally useful. But the \textit{comprehension} of text does not carry with it such pragmatic motivations, nor is it necessarily linear; so the comprehension of logic, a text subset, is not necessarily linear either. There is little about how we represent and comprehend logic that should demand a fundamental one-dimensionality like text beyond its expression through language, and language itself is rarely a systematic and one-to-one parsing of words in sentences to logical semantics. If language were this systematic the reading and writing and parsing of any sentence with a characteristic like Russell's paradox, for example ``I am all of that which I cannot be,'' would be catastrophic to comprehension. Yet despite the paradox of this phrase, we can read it, write it, parse it, and even derive conclusions, albeit absurd ones. Such sentences are simple counterexamples to logic and therefore to proofs with logic having inherently linear comprehension; it is not language that constrains the dimensionality of logic to text, let alone the dimensionality of proof.

\subsection{The shadow of classical logic}

Perhaps it is out of the crooked squiggles of proofs scrawled over folded and sometimes crumpled pages that this inherent one-dimensionality also seems less intrinsic to me. The idea of a ``two-dimensional proof'' is not that alien: proofs by cases, for example, can be seen as ``paths'' of proof satisfied in conjunction, a view that hearkens to Frege's own tree-like notations. But classical logic in its textual symbolic form is a monolith that rules much of our reasoning and so to even suggest re-representation risks backlash. In kind and approaching such risks gently this work is not an attempt to thwart or overturn existing principles, rather to re-cast familiar representations of mathematics to something non-linear, something novel, but more importantly something compatible with Frege's visual conceptions.

\pagebreak

\section{Concept}
\subsection{Non-linear proof representation}
Representing proof in a non-linear sense requires a more general understanding of what a proof really is. Meta-logic has a lot to say here but the core of what we need is the concept of proof as sequences of sentences, which on their own are linear, or linear \textit{parallel} at most in proofs by cases, combined with two ideas: that proofs are capable of sharing a namespace of objects within those sequences of sentences---which is common in bi-conditionals---and that the ordering of sentences within sequences is a relationship between sentences unto itself. The overlapping of proofs then makes for non-linear arrangements through relationships of objects instead of just sentences about objects and these higher-order relationships can be represented as graphs. The graph representation extends us a number of desirable concepts: edge \textit{directionality}, edge \textit{colouring} (or \textit{labeling} if you will) and vertex \textit{adjacency}, which may not seem like much for meta-logic's sake, but it might be enough to represent some mathematics.

\subsection{Re-representing AND}

Consider the entailment of a sentence C by two preceding statements A, B. If only A or B is necessary for C to be true, it is common to eliminate one of A or B. But this sort of elimination is redundant when we only talk about sufficiency: anywhere we can eliminate a disjunction (A or B) it is just as well to say that the truth of (A and B) entail C with redundancy in the truth of one of A or B ``through'' the non-redundancy of the other. Going further it just as well to say that the truth of A \textit{thru} B entails C, or equivalently, the truth of B \textit{thru} A entails C; I am highlighting that redundant thru-like additions such as B to the entailment (A entails C) do not actually change the truth of the entailment but they introduce intermediate structure that we might use. Really the entailment A through B entails C represents C through a causality of the sufficiency of B through the sufficiency of A. If A were sufficient in place of (A or B) then B is sufficient through A, specifically sufficiently redundant through the non-redundancy of A, so C is sufficient through A. With any conjunction (A and B) it is the same, only that we expect the sufficiency of both A and B without redundancy. But this does not change the way we can represent it: it is just as well to re-represent (A and B) entails C as A \textit{thru} B entails C, or symmetrically B \textit{thru} A entails C.

\subsection{e.g. Equivalence of definitions: $G$ is Abelian.}
Graphically re-representing how to prove the equivalence of two definitions for an Abelian group $G$ proceeds as follows:

\includegraphics[width=0.91\columnwidth]{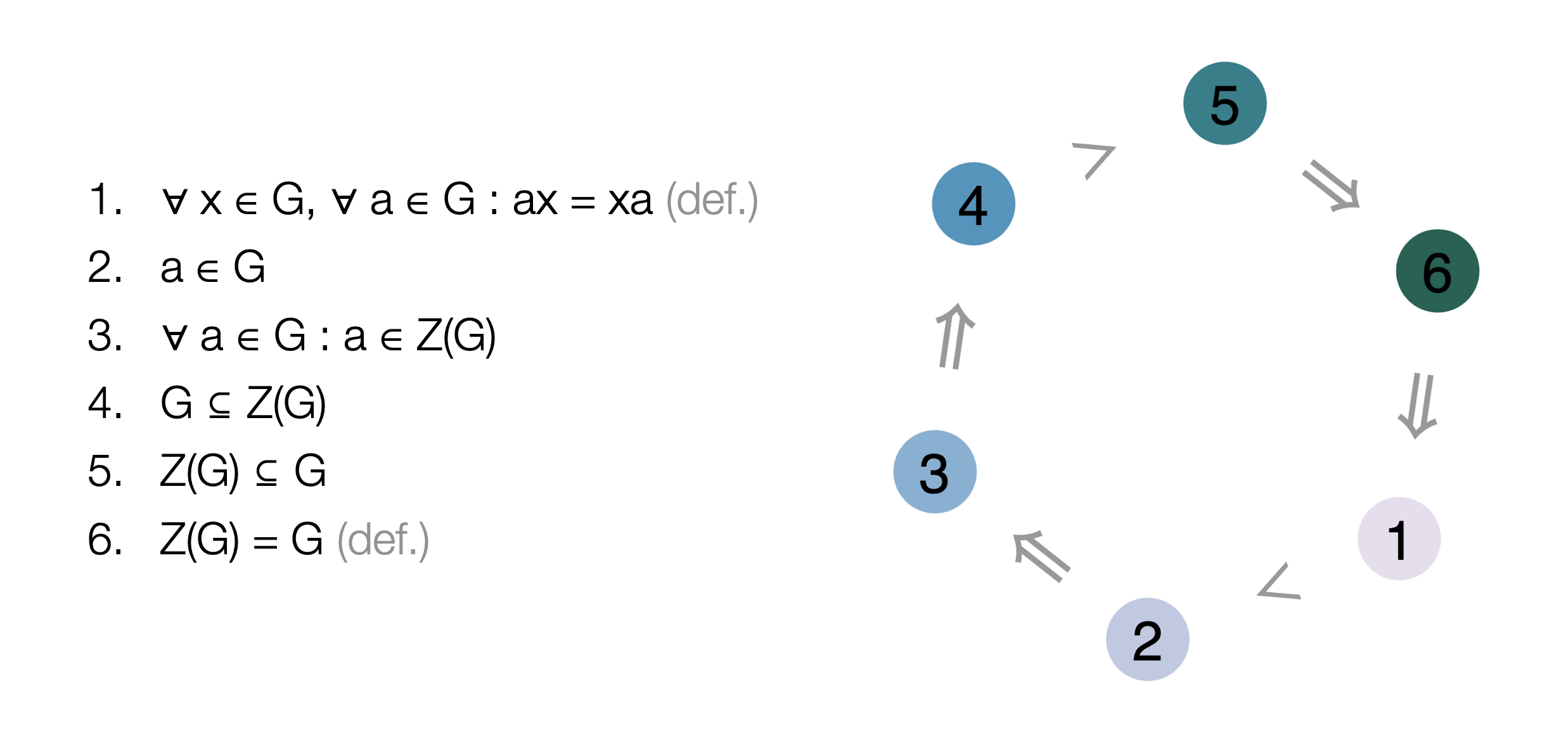}

\subsection{Proof by contraction}
We can speculate from the geometry alone how we might use graph properties to prove if-and-only-if proofs or at least to assert something about their provability. If we overlook edge colours then the uncoloured graph is planar without edge intersections and has a closed circular geometry with no internal edges, so we can contract the entirety of it to a single vertex. This suggests the following rule: 

\begin{quote}
\textit{A set of sentences $\{S\}$ is provable if-and-only-if its graph representation $G=(V,E)$ as vertices $\{V\}$ and edges $\{E\}$ contracts to a single vertex $G_0 = (\{v_0\}, \varnothing\}$.}
\end{quote}

In other words, $\{S\}$ is provable if-and-only-if $G$ is \textit{homomorphic} to $G_0$, or equivalently provable if $G$ is iteratively \textit{contractible} to $G_0$. Maybe this is useful.

\section{Conclusion}
I introduced \textit{Fregean flows}, a graph theoretic representation of classical logic, and showed that it might be useful to prove or perhaps at best show the provability of simple deductive statements. I went further to show that an ``if-and-only-if'' proof under this representation has the geometry of the circle and that we can speculate on graph properties like edge contraction to determine its provability.

\bibliographystyle{plain}
\bibliography{paper}

\end{document}